\magnification 1200
  \input amssym
  \input miniltx
  \input pictex

  %
  \font \bbfive = bbm5
  \font \bbseven = bbm7
  \font \bbten = bbm10
  \font \eightbf = cmbx8
  \font \eighti = cmmi8 \skewchar \eighti = '177
  \font \eightit = cmti8
  \font \eightrm = cmr8
  \font \eightsl = cmsl8
  \font \eightsy = cmsy8 \skewchar \eightsy = '60
  \font \eighttt = cmtt8 \hyphenchar \eighttt = -1

  \font \sixi = cmmi6 \skewchar \sixi = '177
  \font \sixrm = cmr6
  \font \sixsy = cmsy6 \skewchar \sixsy = '60
  \font \tensc = cmcsc10

  \scriptfont \bffam = \bbseven
  \scriptscriptfont \bffam = \bbfive
  \textfont \bffam = \bbten

  \newskip \ttglue

  \def \eightpoint {\def \rm {\fam 0 \eightrm }\relax
  \textfont 0= \eightrm
  \scriptfont 0 = \sixrm \scriptscriptfont 0 = \fiverm
  \textfont 1 = \eighti
  \scriptfont 1 = \sixi \scriptscriptfont 1 = \fivei
  \textfont 2 = \eightsy
  \scriptfont 2 = \sixsy \scriptscriptfont 2 = \fivesy
  \textfont 3 = \tenex
  \scriptfont 3 = \tenex \scriptscriptfont 3 = \tenex
  \def \it {\fam \itfam \eightit }\relax
  \textfont \itfam = \eightit
  \def \sl {\fam \slfam \eightsl }\relax
  \textfont \slfam = \eightsl
  \def \bf {\fam \bffam \eightbf }\relax
  \textfont \bffam = \bbseven
  \scriptfont \bffam = \bbfive
  \scriptscriptfont \bffam = \bbfive
  \def \tt {\fam \ttfam \eighttt }\relax
  \textfont \ttfam = \eighttt
  \tt \ttglue = .5em plus.25em minus.15em
  \normalbaselineskip = 9pt
  \def \MF {{\manual opqr}\-{\manual stuq}}\relax
  \let \sc = \sixrm
  \let \big = \eightbig
  \setbox \strutbox = \hbox {\vrule height7pt depth2pt width0pt}\relax
  \normalbaselines \rm }

  \def \withfont #1#2{\font \auxfont =#1 {\auxfont #2}}

  %

  \def \TRUE {Y}
  \def \FALSE {N}
  \def \EMPTY {}

  \def \ifundef #1{\expandafter \ifx \csname #1\endcsname \relax }

  \def \undefrule{\kern 2pt \vrule width 5pt height 5pt depth 0pt \kern 2pt}
  \def\UndefFlag{}
  \def \possundef #1{\ifundef {#1}\undefrule {\eighttt #1}\undefrule
    \global \edef \UndefFlag{\UndefFlag#1\par }
  \else \csname #1\endcsname \fi }

  %

  \newcount \secno \secno = 0
  \newcount \stno \stno = 0
  \newcount \eqcntr \eqcntr = 0

  \ifundef {showlabel} \global \def \showlabel {\FALSE} \fi  
  \ifundef {auxfile}   \global \def \auxfile   {\TRUE} \fi

  \def \define #1#2{\global \expandafter \edef \csname #1\endcsname {#2}}
  \def \error #1{\parindent 0pt \medskip \bf *******\hfil *******\hfil *******\hfil *******\hfil *******\hfil
    *******\hfil *******\hfil *******\break #1. Exiting... \end }

  \def \advseqnumbering {\global \advance \stno by 1 \global \eqcntr =0}

  \def \current {\ifnum \secno = 0 \number \stno \else \number \secno \ifnum \stno = 0 \else .\number \stno \fi \fi}

  \begingroup \catcode `\@=0 \catcode `\\=11 @global@def@textbackslash{\} @endgroup

  %
  \def \deflabel #1#2{%
    \if\TRUE\showlabel \hbox {\sixrm [[ #1 ]]} \fi
    \ifundef {#1PrimarilyDefined}%
      \define{#1}{#2}%
      \define{#1PrimarilyDefined}{#2}%
      \if\TRUE\auxfile \immediate \write 1 {\textbackslash newlabel {#1}{#2}}\fi
    \else
      \edef \old {\csname #1\endcsname}%
      \edef \new {#2}%
      \if \old \new \else \error{Duplicate definition for label ``{\tt #1}'', already defined as ``{\tt \csname #1\endcsname}''}\fi
      \fi}

  \def \label #1 {\deflabel {#1}{\current }}

  \def \equationmark #1 {\ifundef {InsideBlock}
	  \advseqnumbering
	  \eqno {(\current )}
	  \deflabel {#1}{\current }
	\else
	  \global \advance \eqcntr by 1
	  \edef \subeqmarkaux {\current .\number \eqcntr }
	  \eqno {(\subeqmarkaux )}
	  \deflabel {#1}{\subeqmarkaux }
	\fi }

  \def \split #1.#2.#3.#4;{\global \def \parone {#1}\global \def \partwo {#2}\global \def \parthree {#3}\global \def \parfour {#4}}
  \def \NA {NA}
  \def \ref #1{\split #1.NA.NA.NA;(\possundef {\parone }\ifx \partwo \NA \else .\partwo \fi )}

  %
  \newcount \bibno \bibno = 0

  \def \Bibitem #1 #2; #3; #4 \par{\smallbreak
    \global \advance \bibno by 1
    \item {[\possundef{#1}]} #2, {``#3''}, {#4}.\par
    \ifundef {#1PrimarilyDefined}\else
      \error{Duplicate definition for bibliography item ``{\tt #1}'', already defined in ``{\tt [\csname #1\endcsname]}''}
      \fi
	\ifundef {#1}\else
	  \edef \prevNum{\csname #1\endcsname}
	  \ifnum \bibno=\prevNum \else
		\error{Mismatch bibliography item ``{\tt #1}'', defined earlier as ``{\tt \prevNum}'' but should be ``{\tt \number\bibno}''}
		\fi
	  \fi
    \define{#1PrimarilyDefined}{#2}%
    \if\TRUE\auxfile \immediate\write 1 {\textbackslash newbib {#1}{\number\bibno}}\fi}

  \def \jrn #1, #2 (#3), #4-#5;{\sl #1, \bf #2 \rm (#3), #4--#5}
  \def \Article #1 #2; #3; #4 \par{\Bibitem #1 #2; #3; \jrn #4; \par}

  \def \references {\begingroup \bigbreak \eightpoint \centerline {\tensc References} \nobreak \medskip \frenchspacing }

  %

  \catcode `\@=11
  \def \c@itrk #1{{\bf \possundef {#1}}} 
  \def \c@ite #1{{\rm [\c@itrk{#1}]}}
  \def \sc@ite [#1]#2{{\rm [\c@itrk{#2}\hskip 0.7pt:\hskip 2pt #1]}}
  \def \du@lcite {\if \pe@k [\expandafter \sc@ite \else \expandafter \c@ite \fi }
  \def \cite {\futurelet\pe@k \du@lcite }
  \catcode `\@=12

  %
  \def \Headlines #1#2{\nopagenumbers
    \headline {\ifnum \pageno = 1 \hfil
    \else \ifodd \pageno \tensc \hfil \lcase {#1} \hfil \folio
    \else \tensc \folio \hfil \lcase {#2} \hfil
    \fi \fi }}

  \def \title #1{\medskip\centerline {\withfont {cmbx12}{\ucase{#1}}}}

  \long \def \Quote #1\endQuote {\begingroup \leftskip 35pt \rightskip 35pt
\parindent 17pt \eightpoint #1\par \endgroup }
  \long \def \Abstract #1\endAbstract {\bigskip \Quote \noindent #1\endQuote }
  
  \def \Authors #1{\bigskip \centerline {\tensc #1}}
  \def \Note #1{\footnote {}{\eightpoint #1}}
  \def \Date #1 {\Note {\it Date: #1.}}

  \def \part #1#2{\vfill \eject \null \vskip 0.3\vsize
    \withfont{cmbx10 scaled 1440}{\centerline{PART #1} \vskip 1.5cm \centerline{#2}} \vfill\eject }

  %

  \def \fix {\smallskip \noindent $\blacktriangleright $\kern 12pt}
  \def \iskip {\medskip\noindent}

  \def \ucase #1{\edef \auxvar {\uppercase {#1}}\auxvar }
  \def \lcase #1{\edef \auxvar {\lowercase {#1}}\auxvar }

  \def \section #1 \par {\global \advance \secno by 1 \stno = 0
    %
    \goodbreak \bigbreak
    \noindent {\bf \number \secno .\enspace #1.}
    \nobreak \medskip \noindent }

  \def \state #1 #2\par {\begingroup \def \InsideBlock {} \medbreak \noindent \advseqnumbering {\bf \current .\enspace
#1.\enspace \sl #2\par }\medbreak \endgroup }

  \def \definition #1\par {\state Definition \rm #1\par }

  \long \def \Proof #1\endProof {\begingroup \def \InsideBlock {} \medbreak \noindent {\it Proof.\enspace }#1
\ifmmode \eqno \endproofmarker $$ \else \hfill $\endproofmarker $ \looseness = -1 \fi \medbreak \endgroup }

  \def \$#1{#1 $$$$ #1}
  \def \explica #1#2{\mathrel {\buildrel \hbox {\sixrm #2} \over #1}}
  \def \explain #1#2{\explica{#1}{\ref{#2}}}  
  \def \=#1{\explain {=}{#1}}

  \newcount \fnctr \fnctr = 0
  \def \fn #1{\global \advance \fnctr by 1
    \edef \footnumb {$^{\number \fnctr }$}%
    \footnote {\footnumb }{\eightpoint #1\par \vskip -10pt}}

  \def \text #1{\hbox {#1}}

  %
  
  \def \Item #1{\smallskip \item {{\rm #1}}}
  \newcount \zitemno \zitemno = 0

  \def \izitem {\global \zitemno = 0}
  \def \zitemplus {\global \advance \zitemno by 1 \relax }
  \def \rzitem {\romannumeral \zitemno }
  \def \rzitemplus {\zitemplus \rzitem } 
  \def \zitem {\Item {{\rm (\rzitemplus )}}}
  
  \def \zitemmark #1 {\deflabel {#1}{\rzitem }}

  \newcount \nitemno \nitemno = 0
  \def \initem {\nitemno = 0}
  \def \nitem {\global \advance \nitemno by 1 \Item {{\rm (\number \nitemno )}}}

  \newcount \aitemno \aitemno = -1
  \def \boxlet #1{\hbox to 6.5pt{\hfill #1\hfill }}
  \def \iaitem {\aitemno = -1}
  \def \aitemconv {\ifcase \aitemno a\or b\or c\or d\or e\or f\or g\or
h\or i\or j\or k\or l\or m\or n\or o\or p\or q\or r\or s\or t\or u\or
v\or w\or x\or y\or z\else zzz\fi }
  \def \aitem {\global \advance \aitemno by 1\Item {(\boxlet \aitemconv )}}
  \def \aitemmark #1 {\deflabel {#1}{\aitemconv }}

  \def\Bitem{\Item{$\bullet$}}

  \def \Case #1:{\medskip \noindent {\tensc Case #1:}}

  %
  \def \<{\left \langle \vrule width 0pt depth 0pt height 8pt }
  \def \>{\right \rangle }
  \def \({\big (}
  \def \){\big )}
  
  \def \and {\hbox {,\quad and \quad }}
  \def \calcat #1{\,{\vrule height8pt depth4pt}_{\,#1}}
  
  \def \IFF {\kern 7pt\Leftrightarrow \kern 7pt}
  \def \IMPLY {\kern 7pt \Rightarrow \kern 7pt}
  \def \for #1{,\quad \forall \,#1}
  \def \endproofmarker {\square } 
  \def \"#1{{\it #1}\/} \def \umlaut #1{{\accent "7F #1}}
  \def \inv {^{-1}}
  \def \*{\otimes }
  \def \caldef #1{\global \expandafter \edef \csname #1\endcsname {{\cal #1}}}
  \def \bfdef #1{\global \expandafter \edef \csname #1\endcsname {{\bf #1}}}
  \bfdef N \bfdef Z \bfdef C \bfdef R

  %

  \if\TRUE\auxfile
    \IfFileExists {\jobname.bib}{\input \jobname.bib }{\null}
    \IfFileExists {\jobname.aux}{\input \jobname.aux }{\null}
    \immediate \openout 1 \jobname.aux
    \fi

  \def \close {\if \TRUE \auxfile \closeout 1 \fi
    \if \EMPTY \UndefFlag \else
      \message {*** There were undefined labels ***} \iskip
      ****************** \ Undefined Labels: \tt \par \UndefFlag
      \fi
    \par \vfill \supereject \end }

  %

  \def \Caixa #1{\setbox 1=\hbox {$#1$\kern 1pt}\global \edef \tamcaixa {\the \wd 1}\box 1}
  \def \caixa #1{\hbox to \tamcaixa {$#1$\hfil }}

  \def \med #1{\mathop {\textstyle #1}\limits }

  \def \medcup {\med \bigcup }

  %

  \def \src {d}	                
  \def \ran {r}	                
  \def \E {{\cal E}}                               


  \def \acite [#1]{\cite [#1]{actions}}



  \def \S {{\cal S}}
  \def \E {{\cal E}}
  \def \interior #1{{\buildrel \circ \over #1}}

  \def \TF {T\!F}
  \def \J {{\cal J}}
  \def \tight {_{\rm tight}}
  \def \red {_{\rm red}}


  \def \G {{\cal G}}
  \def \Dm #1#2{D^#1_{#2}}
  \def \D #1#2{D^{#1}\{#2\}}  
  \def \act {\alpha}
  \def \actfil {\beta}
  \def \X #1#2{{\cal F}^{\,#1}_{#2}} 
  \def \Gth {\G \tight (\S )}
  \def \Eth {\hat \E \tight }
  \def \zero {^{(0)}}



  \def \orb #1{\hbox {Orb}(#1)}


  \def \clos #1{\overline {#1}}






  \def \mathcal #1{{\cal #1}}
  \def \text #1{\hbox {#1}}
  \def \mathbb #1{{\bf #1}}
  \def \mbox #1{\hbox {#1}}
  \def \frac #1#2{{#1 \over #2}}

  \def \vspace #1{\vskip #1}
  \def \hspace #1{\hskip #1}

  \newcount \itemdpt \itemdpt = 0
  \def \uplevel {\begingroup \global \advance \itemdpt by 1 \advance \parindent by 18pt
  \ifcase \itemdpt \or \initem \or \iaitem \or \izitem \fi }
  \def \dnlevel {\par \global \advance \itemdpt by -1 \advance \parindent by -18pt \endgroup }
  \def \itm {\ifcase \itemdpt \or \nitem \or \aitem \or \zitem \fi }






  \def \quoapprox {E^0{\kern -1pt/\kern -2pt}\approx }

  %

  \makeatletter
  \def \Gin@driver{pdftex.def}
  \input color.sty
  \resetatcatcode

	\def \beginRuyComment #1 \endRuyComment {
	  \goodbreak \begingroup \rm \color {magenta}\bigskip \hrule \parindent 0pt\parskip 10pt COMMENT:\par #1
	  \par \bigskip \hrule \bigskip \endgroup \noindent }



  \font\fixed = cmitt10 
  \long\def\prelim #1\endprelim{{\overfullrule=0cm \fixed (This is preliminary) #1\par}}


\font \titlefont = cmbx10
\def\title#1{\medskip\centerline {\titlefont \ucase{#1}}}

\null\bigskip

\title{THE TIGHT GROUPOID OF AN INVERSE SEMIGROUP}

  \Headlines {THE TIGHT GROUPOID OF AN INVERSE SEMIGROUP}
  {R.~Exel and E.~Pardo}

  \Authors {Ruy Exel and Enrique Pardo}

  \Date {21 August 2014} 


  \Note {\it Key words and phrases: \rm Inverse semigroup, semi-lattice, tight character, tight filter, ultra-filter, groupoid, groupoid C*-algebra.}

  \Note {The first-named author was partially supported by CNPq. The second-named author was partially supported by PAI
III grants FQM-298 and P11-FQM-7156 of the Junta de Andaluc\'{\i}a and by the DGI-MICINN and European Regional
Development Fund, jointly, through Project MTM2011-28992-C02-02.}

\bigskip

  \Abstract In this work we present algebraic conditions on an inverse semigroup $\S$ (with zero) which imply that its
associated
tight groupoid $\Gth$ is: Hausdorff, essentially principal, minimal and contracting, respectively.
  In some cases these conditions  are in fact
necessary and sufficient.
  \endAbstract

\bigskip

\section Introduction

This article should  be
considered as a continuation of the work started  by the first named author in \cite {actions}, where the notion of
tight representations of inverse semigroups was introduced and used to study a large class of C*-algebras.
Given any inverse semigroup  $\S$, the theory developed in \cite{actions} gives a recipe to build an
\'etale groupoid, denoted  $\Gth$, whose C*-algebra is isomorphic to the universal C*-algebra for tight representations of $\S$.

Since then a very large number of C*-algebras were shown to fit this model, including virtually  all
C*-algebras in the literature defined in terms of generators and relations, provided the relations specify that the
generators are partial isometries\fn
  {A relevant exception is the universal C*-algebra generated by a single partial isometry.}.
  To be honest, algebras which have been referred to as \"{Toeplitz extensions} in various contexts are typically not
included, but these  are often extensions of some other algebra of interest which, most of the time, may be
shown to fit the model referred to above.

To mention an example of interest to us,  in  \cite{EP} and \cite{EPTwo} a
unified treatment was given to a certain class of C*-algebras studied by Katsura in \cite{KatsuraOne}, alongside
Nekrashevych's C*-algebras introduced and discussed in \cite{NekraJO}, \cite{Nmsn} and \cite{NC}.  The unifying
principle is the notion of self similar graphs introduced in \cite{EPTwo}, which gives rise to C*-algebras that can be
effectively studied via the already mentioned theory of tight representations of   inverse semigroups.

The applications of this circle of ideas to a growing class of C*-algebras begs for a unified treatment  of questions of
relevance in the study of the structure of such C*-algebras.  Since these are groupoid C*-algebras for $\Gth$, and since
many algebraic properties of the C*-algebra depend directly on this groupoid, we have spent a lot of effort in trying to
characterize  these properties in terms of  the algebraic structure of the given inverse semigroup.
  The present  paper is thus our account of that effort, in which we were more or less successful in determining conditions on $\S$, often
necessary and sufficient, for $\Gth$ to be:
  \Bitem Hausdorff,
  \Bitem essentially principal,
  \Bitem minimal and
  \Bitem contracting.
  \iskip These conditions are all algebraic in nature and their verification is often easily determined, as illustrated
by the applications given  in \cite{EPThree}.

After this paper first appeared in the arXiv, we had access to a preprint by Steinberg \cite{SteinbergPreprint} where
similar questions are treated.

\section Inverse semigroups

The specific purpose of the present section is to review some of the main concepts from \cite {actions}, occasionally
offering minor improvements.  We will also  briefly recall some  basics facts about inverse semigroups.
We refer the reader to \cite {Lawson} for an extensive treatment of inverse semigroups.

\definition \label DefineInvSGroup
  An \"{inverse semigroup} is a set $\S $ equipped with an associative multiplication operation
  $$
  \S \times\S \to\S ,
  $$
  such that:
  \izitem
  \zitem for every $s$ in $\S $, there exists a unique $s^*\in\S $, such that $ss^*s =s$, and $s^*ss^* = s^*$,
  \zitem there exists a (necessarily unique) element $0\in\S $, called the \"{zero element}, such that $s0=0s=0$, for all
$s$ in $\S $.

\medskip Axiom \ref{DefineInvSGroup.ii} is not usually part of the standard definitions of inverse semigroups, but
in all of our uses of this concept, the zero element will play an important role, hence our insistence in including it
explicitly.

In case an inverse semigroup $\S $ lacks a zero, it is always easy to adjoin one by simply considering $\S \cup\{0\}$ with
the obvious extension of the multiplication operation.  In parti\-cu\-lar it is imperative to do so in case one wishes to
consider groups as special cases of inverse semigroups.  Although this might seem a little strange, it causes no serious
technical problems.

Whenever an inverse semigroup $\S $ is in sight, we will denote by
  $$
  \E = \{e\in \S : e^2=e\}.
  $$
  In the rare occasions when explicit reference to $\S $ is important to  dispel possible confusions, we might also write $\E (\S )$ for
the above set.

  It is well known that every element $e$ in $\E $ is \"{self-adjoint} in the sense that $e^*=e$.  Moreover $\E $ is a
commutative sub-semigroup of $\S $ and, under the order relation
  $$
  e\leq f \IFF e = ef \for e,f\in\E ,
  \equationmark OrderForIdempotents
  $$
  one has that the greatest lower bound $e\wedge f$ always exists for any $e$ and $f$ in $\E $, namely
  $$
  e\wedge f=ef.
  $$
  For this reason, $\E $ is called the \"{idempotent semi-lattice} of $\S $.

The order relation on $\E $ may in fact be extended to an order relation on $\S $, defined by
  $$
  s\leq t \IFF s = ts^*s \for s,t\in\S .
  \equationmark DefineOrder
  $$
  It is well known that this relation is invariant under left or right multiplication and that $s\leq t$, if and only if
  $s = ss^*t$.  Moreover,
  $$
  s\leq t \IMPLY s^*s\leq t^*t, \hbox { and } ss^*\leq tt^*.
  \equationmark CompareIniFin
  $$
  See \cite [Proposition 7]{Lawson} for proofs of these facts.

  Given two elements $e$ and $f$ in $\E $, we will say that $e$ is \"{orthogonal} to $f$, provided $ef=0$.  In symbols
  $$
  e\perp f \IFF ef=0.
  $$
  On the other hand, we say that $e$ \"{intersects} $f$, when $ef\neq0$.  In symbols
  $$
  e\Cap f \IFF ef\neq0.
  $$

A \"{character} on $\E $ is any nonzero map
  $$
  \phi:\E \to \{0,1\}
  $$
  such that $\phi(0)=0$, and
  $$
  \phi(ef) = \phi(e) \phi(f),
  $$
  for all $e,f\in\E $.

  The set of all characters on $\E $ is denoted by $\hat \E _0$.  This notation is meant to avoid confusion with the set
$\hat \E $ of all \"{semi-characters}, meaning characters which are not required to satisfy ``$\phi(0)=0$''.  Semi-characters
are used by some authors but they will not play any role in this work.

  We will always view $\hat \E _0$ as a topological space equipped with the product topology, that is, the subspace
topology inherited from $\{0,1\}^\E $.

Should the zero map be allowed as a character, $\hat \E _0$ would be closed in $\{0,1\}^\E $, and hence compact.  Having
explicitly excluded the zero map from $\hat \E _0$, this becomes a locally compact space.

A \"{filter} in $\E $ is a nonempty subset $\eta\subseteq\E $, not containing the zero element, which is closed under ``$\wedge $'',
and which moreover satisfies
  $$
  f\geq e\in\eta \IMPLY f\in\eta,
  $$
  for all $e,f\in\E $.
  Given a filter $\eta$, we define
  $$
  \phi_\eta : e \in \E \mapsto [e\in\eta]\in\{0,1\},
  $$
  where the brackets correspond to boolean value.  In other words, $\phi_\eta$ is the characteristic function of $\eta$, when $\eta$
is seen as a subset of $\E $.  It is easy to see that $\phi_\eta$ is a character on $\E $.

  Conversely, given a character $\phi$ on $\E $, the set
  $$
  \eta_\phi = \{e\in\E : \phi(e)=1\}
  $$
  is a filter.  One may easily prove that the correspondences ``$\phi\to\eta_\phi$'' and ``$\eta\to\phi_\eta$'' are each other's inverse, and
hence filters may, and will, be identified with characters and vice versa.

The product topology viewed from the point of view of filters is easy to describe: given finite subsets $X,Y\subseteq\E $,
consider the set
  $$
  U(X,Y) = \{\eta\in\hat \E _0: X\subseteq\eta,\ Y\subseteq\E \setminus \eta\}.
  $$
  Then each $U(X,Y)$ is an open set and the collection of all such is easily seen to form a basis for the topology of
$\hat \E _0$.  Assuming that $X$ is nonempty, and letting
  $x_0 = \bigwedge X$, observe that for any filter $\eta$ one has that
  $$
  x_0\in\eta \iff X\subseteq\eta,
  $$
  so
  $$
  U(X,Y) = U\big(\{x_0\},Y\big).
  $$

  On the other hand, if $X$ is empty, and if $\eta$ is a given element in $U(\emptyset ,Y)$ then, choosing any $e\in\eta$, one
has that
  $$
  \eta \in U\big(\{e\},Y\big) \subseteq U(\emptyset ,Y).
  $$
  This shows that
  $$
  U\big(\emptyset ,Y\big) = \medcup _{e\in\E } U\big(\{e\},Y\big),
  $$
  and we then see that the collection of all $U(X,Y)$, where $X$ is a singleton, also form a basis for the topology of
$\hat \E _0$.

A filter $\xi$ is said to be an \"{ultra-filter} if it is not properly contained in another filter.  A character $\phi$ is
said to be an \"{ultra-character} if its associated filter $\xi_\phi$ is an ultra-filter.

A useful characterization of ultra-filters is given by \acite [Lemma 12.3]: a filter $\xi$ is an ultra-filter if and only
if $\xi$ contains every idempotent $f$ such that $f\Cap e$ for every $e$ in $\xi$.  Therefore the only reason why an
idempotent $e$ fails to belong to an ultra-filter $\xi$ is when $e$ is orthogonal to some $f$ in $\xi$.

Referring to our discussion above regarding open subsets of $\hat \E _0$, ultra-filters have specially nice neighborhood
bases:

\state Proposition \label NicerNeighborhood
  If $\xi$ is an ultra-filter then the open sets of the form
  $
  U\big(\{e\},\emptyset \big),
  $
  as $e$ range in $\xi$,
  forms a neighborhood basis for $\xi$.

\Proof
  It is enough to show that, whenever $X$ and $Y$ are finite subsets of $\E $ such that $\xi\in U(X,Y)$, then there is some
$e$ in $\xi$ such that
  $$
  \xi\in U\big(\{e\},\emptyset \big) \subseteq U(X,Y).
  \equationmark GoalSimpleNbdOfUltre
  $$

  For each $y$ in $Y$, we have that $y\notin\xi$ so, by \acite [Lemma 12.3], there is some $f_y\in\xi$ such that $f_y\perp y$.  Defining
  $$
  e = \bigwedge X \wedge \bigwedge _{y\in Y}f_y,
  $$
  the reader may now easily show that $e$ satisfies \ref{GoalSimpleNbdOfUltre}, and hence the proof is concluded.
\endProof

\definition \label DefineIdeal
  A subset $\J \subseteq\E $ is said to be an \"{ideal} of $\E $ if
  \izitem
  \zitem $0\in\J $, and
  \zitem for every $e$ in $\J $, and every $f$ in $\E $, one has that $ef\in\J $.

Notice that \ref{DefineIdeal.ii} is equivalent to saying that
  $$
  f\leq e \in \J \IMPLY f\in\J ,
  $$
  for every $e,f\in\E $.
  Ideals are therefore precisely the nonempty hereditary subsets of $\E $.

Important examples are the \"{principal ideals}, namely ideals of the form
  $$
  \J _e := \{f\in\E : f\leq e\} = e\E .
  \equationmark PrincipalIdeal
  $$

  Like the similar concept in ring theory, an arbitrary intersection of ideals is an ideal, but, unlike rings, an
arbitrary union of ideals in a semi-lattice is also an ideal.

If $\J $ is an ideal in $\E $, then so is
  $$
  \J ^\perp = \{f\in\E : f\perp e, \hbox { for all } e\in\J \}.
  $$
  Given any $e$ in $\E $, notice that
  $$
  \J _e^\perp = \{f\in\E : f\perp e\}.
  $$

  If $X$ and $Y$ are subsets of $\E $, one then has that
  $$
  \E ^{X,Y} := \med\bigcap _{x\in X} \J _x \cap \med\bigcap _{y\in Y} \J _y^\perp
  \equationmark DefineEXY
  $$
  is an ideal of $\E $, which will soon play an important role in the definition of a special property for characters.
Notice that if $X$ is finite and nonempty, and if we let $x_0 = \bigwedge X$, then
  $$
  \E ^{X,Y} = \E ^{\{x_0\},Y}.
  $$

\definition \label DefineCover
  Given an ideal $\J \subseteq\E $, and a subset $C\subseteq\E$,  we will  say that:
  \iaitem
  \aitem $C$ is an \"{outer cover} for $\J $ if, for every nonzero $f$ in $\J $, there exists some $c$ in $C$ such that
$c\Cap f$,
  \aitem $C$ is a \"{cover} for $\J $ if $C$ is an outer cover for $\J $ and $C\subseteq\J$.
  \iskip
  Given an idempotent $e$ in $\E$, and a cover (resp.~outer cover) $C$ for $\J _e$, we will say that $C$ is a
\"{cover}
(resp.~\"{outer cover}) for $e$.

Covers were introduced in \acite [Definition 11.5], where the subset $\J $ being \"{covered} was not required to be an
ideal.  However only covers for ideals will be relevant here.

It is easy to see that $C$ is an outer cover for an idempotent $e$ if and only if
  $$
  Ce:=\{ce: c\in C\}
  $$
  is a cover for $e$.  However it does not seem possible to turn an outer cover of a general ideal into a
cover, especially if the ideal is of the form $\J _e^\perp $.

A character $\xi$ is said to be a \"{tight character} if, for every finite subsets $X,Y\subseteq\E $, and for every finite cover
  $Z\subseteq \E ^{X,Y}$, one has that
  $$
  \med \bigvee _{z\in Z} \phi(z) \geq \med \bigwedge _{x\in X} \phi(x) \wedge \med \bigwedge _{y\in Y} \big(1-\phi(y)\big).
  $$

Replacing ``$\geq$'' by ``$\leq$'' above, it is elementary to check that the resulting inequality is always true.  When $\phi$ is
tight one therefore gets an equality above.

  We will say that a filter $\xi$ is a \"{tight filter} when its associated character $\phi_\xi$ is a tight character.
Thus, a filter $\xi$ is a tight filter if and only if, for every finite subsets $X,Y\subseteq\E $, and for every finite cover
  $Z\subseteq \E ^{X,Y}$, one has that
  $$
  X\subseteq\xi, \hbox { and } Y\cap\xi=\emptyset  \IMPLY   Z\cap\xi\neq\emptyset .
  \equationmark TightnessForFilters
  $$

The set of all tight characters is called the \"{tight spectrum} of $\E $, and is denoted by $\hat \E \tight $.  It is
easy to see that $\hat \E \tight $ is closed within $\hat \E _0$, and hence a locally compact space.

Every ultra-character is known to be tight \acite [Proposition 12.7], and the set of all ultra-characters is dense in
$\hat \E \tight $ \acite [Theorem 12.9].

\section Inverse semigroup actions

In this section we will focus on actions of inverse semigroups on locally compact topological spaces.  Besides
introducing the basic concepts, we will briefly recall the definition of the groupoid of germs for a given inverse
semigroup action.  Considering that such a groupoid is often non-Hausdorff, we will discuss conditions under which this
pathology does not occur.  We will then give necessary and sufficient conditions for the tight groupoid associated to an
inverse semigroup to be Hausdorff.  We again refer the reader to \cite {actions} for the basic theory of inverse
semigroup actions.

\definition \label DefineAction
  Throughout this work we will write
  $$
  \act :\S \curvearrowright X
  $$
  to mean that
  \izitem
  \zitem $\S $ is an inverse semigroup (with zero),
  \zitem $X$ is a locally compact Hausdorff topological space,
  \zitem $\act $ is an action of $\S $ on $X$, in the sense of \acite [Definition 4.3],
  \zitem $\act _0$ is the empty map.

For the benefit of the reader we recall that \ref{DefineAction.iii} means that for each $s\in\S $, one is given a
\"{partial homeomorphism} $\act _s$ on $X$ (meaning a homeomorphism between two open subsets of $X$), such that
  $$
  \act _s\circ\act _t = \act _{st} \for s,t\in\S ,
  $$
  where the composition above is defined on the largest domain where it makes sense.

It follows that $\act _{s^*} = \act _s\inv $, for every $s$ in $\S $.  Moreover, for every $e$ in $\E $, the associated
partial homeomorphism $\act _e$ is necessarily the identity map on some open subset of $X$, which we will always denote
by
  $\Dm \act e$.  The collection of sets
  $$
  \{\Dm \act e\}_{e\in\E }
  \equationmark IntroduceDomains
  $$
  should therefore be thought of as an important ingredient of $\act $.  Incidentally, when describing the definition of
an action \acite [Definition 4.3] above, we should have added that the $\Dm \act e$ are required to cover $X$.

  It may then be easily verified that, for every $s$ in $\S$, the domain of  $\act _s$ coincides with $\Dm \act {s^*s}$, while the range of
$\act _s$ coincides with $\Dm \act {ss^*}$.  In other words,
  $$
  \act _s: \Dm \act {s^*s} \to \Dm \act {ss^*}.
  $$

If $Y$ is a subset of $X$, we say that $Y$ is \"{invariant} under $\act $ if
  $$
  \act _s(Y\cap\Dm \act {s^*s})\subseteq Y \for s \in \S .
  $$

  If $Y$ is moreover  locally compact, we may restrict $\act $ to $Y$, thus obtaining an action of $\S $ on $Y$.

Given any inverse semigroup $\S $, there is an important action
  $$
  \beta: \S \curvearrowright \hat \E _0,
  $$
  which we would like to describe.  First of all, for each $e$ in $\E $, let
  $$
  \Dm \actfil e = \{\phi\in\hat \E _0: \phi(e) =1\}.
  $$

  For each $s$ in $\S $, we let $\actfil _s$ be the map from $\Dm \actfil {s^*s}$ to $\hat \E _0$ given by
  $$
  \actfil _s(\phi)\calcat e = \phi(s^*es) \for \phi\in\Dm \actfil {s^*s} \for e\in \E .
  $$
  Observe that $\actfil _s(\phi)$ is not the zero map because
  $$
  \actfil _s(\phi)\calcat {ss^*} =
  \phi\big(s^*(ss^*)s\big)=
  \phi(s^*s)=1\neq0.
  $$

  It may be shown that the range of $\actfil _s$ is $\Dm \actfil {ss^*}$, and that the correspondence $s\mapsto\actfil _s$ is a
well defined action of $\S $ on $\hat \E _0$.

It is sometimes useful to view $\actfil $ as acting on filters, once these are identified with
characters.  The resulting picture is as follows: for each $e$ in $\E $, one has that
  $$
  \Dm \actfil e = \{\xi\in\hat \E _0: e\in\xi\},
  \equationmark DefineDomFilter
  $$
  while,
  for each $s$ in $\S $, and for each $\xi$ in $\Dm \actfil {s^*s}$, one has that
  $$
  \actfil _s(\xi) = \{f\in\E : f\geq ses^*, \hbox { for some } e\in \xi\}.
  \equationmark DefineActionOfFilters
  $$

\state Proposition \label ActionPreserveUltra
  Given an inverse semigroup $\S $ and an ultra-filter $\xi\in\Dm \actfil {s^*s}$, one has that $\actfil _s(\xi)$ is an
ultra-filter.

\Proof
  Suppose that $\eta$ is a filter containing $\actfil _s(\xi)$.  Then
  $$
  \xi = \actfil _{s^*}\big(\actfil _s(\xi)\big) \subseteq \actfil _{s^*}(\eta),
  $$
  so we conclude that $\xi =\actfil _{s^*}(\eta)$ by maximality, and then
  $$
  \eta = \actfil _s\big(\actfil _{s^*}(\eta)\big) = \actfil _s(\xi),
  $$
  proving that $\actfil _s(\xi)$ is an ultra-filter.  \endProof

Using the above result, or explicitly referring to \acite [Proposition 12.11], we have that $\hat \E \tight $ is also
invariant under $\actfil $.  Restricting $\beta$ to $\hat \E \tight $ we obtain an action
  $$
  \theta:\S \curvearrowright \hat \E \tight
  \equationmark IntroduceTightAction
  $$
  which is by far the most important inverse semigroup action in this work, as it will soon become clear.  We will refer
to it as the \"{standard action} of $\S$.

\state Proposition \label CoversAndCovers
  Let $\S $ be an inverse semigroup with idempotent semi-lattice $\E $.  For  any subset $F\subseteq\E $, let
  $$
  \D \theta F := \medcup _{f\in F}\Dm \theta f.
  $$
  Given an ideal $\J \subseteq\E $, and a finite subset $C\subseteq\E $, one has that   $C$ is an outer cover for $\J $ if and only if
  $$
  \D \theta\J \subseteq \D \theta C .
  $$
  If moreover $C\subseteq\J$, then $C$ is a cover for $\J $ if and only if the above inclusion of sets is an equality.

\Proof Assuming  that $C$ is an outer cover for $\J $,
  take any $\xi$ in $\D \theta\J $, so there is some $e$ in $\J $, such that $\xi\in\Dm \theta e$, and hence $e\in\xi$.  Observe that
  $$
  Ce:= \{ce:c\in C\}
  $$
  is clearly a cover for
  $$
  \J _e = \E ^{\{e\},\emptyset }
  $$
  and, since $\xi$ is tight, we deduce from \ref{TightnessForFilters} that $ce\in\xi$, for some $c$ in $C$.  However, as
$c\geq ce\in\xi$, we have that $c\in\xi$, whence $\xi\in\Dm \theta c\subseteq\D \theta C$.  This proves that $\D \theta\J \subseteq \D \theta C$.

Conversely, let $e$ be a nonzero element in $\J $.  Our task is then to show that there exists some $c$ in $C$, such
that $c\Cap e$.  Using Zorn's Lemma, let $\xi$ be an ultra-filter containing $e$.

 By \acite [Proposition 12.7], we have
that $\xi$ is tight, and then it is clear that $\xi\in\Dm \theta e\subseteq\D \theta\J $.  By hypothesis we then have that $\xi$ is in $\D \theta C$,
which means that $\xi\in\Dm \theta c$, for some $c$ in $C$, whence $c\in\xi$.  Using that $\xi$ is a filter it follows that $ec\in\xi$, so
$ec\neq0$, and then $e\Cap c$, as desired.

  We leave the easy proof of the last sentence in the statement for the reader.
 \endProof

The hypothesis that the subset $C$ is finite in the above result cannot be removed, as one may easily show through
elementary examples.  For that reason only finite covers are considered in this work.  In fact, the very use of the word
``cover'' in the infinite case might not be entirely appropriate and perhaps one should even require finiteness in
Definition \ref{DefineCover}.  In any case only finite covers will play a role  in the sequel.

For the special case of the ideal $\J _e$, notice that $\D \theta{\J _e} = \Dm \theta e$, so Proposition \ref{CoversAndCovers}
takes the following somewhat simpler form:

\state Proposition \label CoversAndCoversForIdempotents
  Given $e\in\E $, and a finite subset $C\subseteq\E$, one has that
  $C$ is an outer cover for $e$ if and only if
  $$
  \Dm \theta e \subseteq \medcup _{c\in C}\Dm \theta c.
  $$
  If moreover $C\subseteq\J_e$, then $C$ is a cover for $\J_e $ if and only if the above inclusion of sets is an equality.

Given any action $\act :\S \curvearrowright X$, let us now briefly describe its groupoid of germs.  The reader is
referred to \cite [page 140]{pat} and \acite [Section 4] for more details.

We begin by considering the set
  $$
  \Omega = \{(s,x)\in\S \times X: x\in\Dm \act {s^*s}\}.
  $$
  Given $(s,x)$ and $(t,y)$ in $\Omega$, we say that
  $$
  (s,x)\sim (t,y)
  $$
  provided $x=y$, and there exists an idempotent $e\in\E $, such that $x\in\Dm \act e$, and $se=te$.  In this case we say
that $s$ and $t$ have the \"{same germ} at $x$.  En passant notice that one then has $\act _s(x) = \act _t(x)$,
because
  $$
  \act _s(x) = \act _s\big(\act _e(x)\big) = \act _{se}(x) = \act _{te}(x) = \act _t\big(\act _e(x)\big) = \act _t(x).
  $$

It is elementary to check that ``$\sim $'' is an equivalent relation on $\Omega$.  The equivalence class of each $(s,x)$ in
$\Omega$, usually denoted by $[s,x]$,  is called the \"{germ} of $s$ at $x$.  The set of all germs, namely
  $$
  \G _\act = \Omega/\sim ,
  $$
  is the carrier set of the groupoid we are about to define.  First of all, given any $[s,x]$ in $\G _\act $, we define
the \"{source}, or \"{domain} of $[s,x]$, as well as its \"{range}, by
  $$
  \src ([s,x]) = x \and \ran ([s,x]) = \act _s(x),
  $$
  respectively.

When $[s,z]$ and $[t,x]$ are given in $\G _\act $, and $\src ([s,z]) = \ran ([t,x])$, that is, when $z = \act _t(x)$, we
define their \"{product}, by
  $$
  [s,z][t,x] = [st,x].
  $$

It is also customary to denote a germ $[t,x]$ by the alternative notation $[y,t,x]$, where $y=\act _t(x)$.  The problem
with this notation is that it carries a superfluous information in the variable $y$, but the definition of the product
looks a bit more symmetric:
  $[w,s,z][y,t,x]$ is defined if and only if $z=y$, in which case
  $$
  [w,s,z][y,t,x]= [w,st,x].
  $$

We leave it up to the reader to decide which is the nicer notation, but in the meantime we shall adopt the shorter
notation $[t,x]$.

It may then be proved that $\G _\act $ is a groupoid with the above operation, the
  inverse\fn
  {The inverse also looks nicer in the alternative notation: $[y,t,x]\inv =[x,t^*,y]$.}
  of any given germ $[s,x]$ being given by
  $$
  [s,x]\inv = [s^*,\act _s(x)].
  $$

  The unit space of $\G _\act $ is given by
  $$
  \G _\act \zero = \big \{[e,x]: e\in\E ,\ x\in\Dm \act e\big \},
  $$
  and it may be identified with $X$ via the bijective mapping
  $$
  [e,x] \in \G _\act \zero \mapsto x \in X.
  \equationmark UnitSpace
  $$

  In order to describe the appropriate topology on $\G _\act $, we introduce the following notation: given $s$ in $\S $,
and any open set $U\subseteq\Dm \act {s^*s}$, we let
  $$
  \Theta(s,U) = \big \{[s,x]: x\in U\big \}.
  $$

  The collection of all $\Theta(s,U)$ may then be shown to form a basis for a topology on $\G _\act $, with respect to which
$\G _\act $ becomes a locally compact, \'etale groupoid.

Even if $X$ is assumed to be Hausdorff, $\G _\act $ is not always Hausdorff.  The question of whether or not
$\G _\act $ is Hausdorff is thus the first main problem we wish to analyze.  We begin by presenting a proof of a well
known characterization of the Hausdorff property for general \'etale groupoids.

\state Proposition \label ThmHausdorffGroupoid
  Let $\G $ be an \'etale groupoid with range and source maps denoted by $\ran $ and $\src $, respectively, and such that
$\G ^{(0)}$ is Hausdorff. Then, the following are equivalent:
  \initem
  \nitem $\G $ is Hausdorff.
  \nitem $\G ^{(0)}$ is closed.

\Proof Assuming that $\G $ is Hausdorff, observe that
  $$
  \G ^{(0)}=\{\gamma\in\G : \gamma=r(\gamma)\}.
  $$
  Therefore $\G ^{(0)}$ is the inverse image, under the continuous mapping
  $$
  f: \gamma\in\G \mapsto (\gamma, r(\gamma)) \in \G \times\G
  $$
  of the diagonal of $\G $.  Since the diagonal is closed in a Hausdorff space, we deduce that $\G ^{(0)}$ is closed.

Conversely, suppose that $\G ^{(0)}$ is closed.  Given two distinct elements $\gamma$ and $\gamma'$ in $\G $, suppose first that $\ran
(\gamma)\neq\ran (\gamma')$.  Then, since $\G ^{(0)}$ is Hausdorff, we may find disjoint open subsets of $\G ^{(0)}$ separating $\ran (\gamma)$ from
$\ran (\gamma')$, and then their inverse image under the continuous function $\ran $ will separate $\gamma$ from $\gamma'$.

Suppose now that $\ran (\gamma)=\ran (\gamma')$.  Then the product $\gamma\inv \gamma'$ is well defined, and since $\gamma\neq\gamma'$, we have that
$\gamma\inv \gamma'\not \in\G ^{(0)}$.  By hypothesis, there is an open neighborhood $V$ of $\gamma\inv \gamma'$ which does not intersect $\G ^{(0)}$
and, because $\G $ is \'etale, we may assume that $V$ is an open bisection. Given any open bisection $W$ such that $\gamma\in W$,
we have that $WV$ is open, and
  $$
  \gamma' = \gamma (\gamma\inv \gamma') \in WV.
  $$

  The proof will then be complete once we show that $W\cap WV=\emptyset $.  To settle this, suppose by contradiction that
$\eta\in W\cap WV$.  We may then write $\eta=\omega\nu $, for suitable $\omega\in W$ and $\nu \in V$, so
  $r(\eta)=r(\omega)$.  But since both $\eta$ and $\omega$ lie in the bisection $W$, where the range map is one-to-one, we deduce that
$\eta=\omega$, whence
  $$
  \nu = \omega\inv \eta = \omega\inv \omega \in\G ^{(0)},
  $$
  and then $\nu \in V\cap\G ^{(0)}$, which is a contradiction.
  \endProof

As already announced, our first major goal is to study the Hausdorff property for the groupoid of germs $\G _\act $
relative to an inverse semigroup action $\act $.  Some important tools for this task are described next.

\definition \label DefineXs
  Given an inverse semigroup $\S $, and given $s$ in $\S $, we will let
  $$
  \J _s = \{e\in\E : e\leq s\}.
  $$
  If we are moreover given an action $\act :\S \curvearrowright X$, we will let
  $$
  \X \act s = \medcup _{e\in\J _s}\Dm \act e.
  $$

Incidentally, regarding the relation ``$e\leq s$'' appearing in the definition of $\J _s$, above, notice that, when comparing
an idempotent element $e$ to a general element $s$ in $\S $, we have by \ref{DefineOrder} that
  $$
  e\leq s \IFF e = se.
  \equationmark compareIdempWithNonIdemp
  $$

It is easy to see that  $\J _s$ is an ideal of $\E $.  However, since $s$ does not necessarily belong to $\E $, we cannot say that
$\J _s$ is a principal ideal as in \ref{PrincipalIdeal}.

Given an action $\act :\S \curvearrowright X$, using the terminology introduced in \ref{CoversAndCovers}, observe
that
  $$
  \X \act s = \D \act {\J _s}.
  \equationmark XsisDJs
  $$
  Also notice that if $t\leq s$, then necessarily $\Dm \act {t^*t}\subseteq\Dm \act {s^*s}$, so we see that
  $$
  \X \act s\subseteq\Dm \act {s^*s} \for s\in \S .
  $$

It should also be noted that if $s$ is idempotent, then the above inclusion of sets becomes an equality.
  On the opposite extreme, if $s$ is not idempotent, and if $\S $ is
  E*-unitary\fn
  {An inverse semigroup $\S $ is said to be \"{E*-unitary} if, whenever an element $s$ in $\S $ dominates a nonzero
idempotent, then $s$ must itself be idempotent.  An alternative way to express this condition is that $\J _s=\{0\}$,
whenever $s$ is not idempotent.},
  then $\J _s=\{0\}$, whence $\X \act s=\emptyset $.

\state Proposition \label TrivialGermsInSlice
  Given an action $\act :\S \curvearrowright X$,
  for every $s$ in $\S $, one has that
  $$
  \Theta(s,\Dm \act {s^*s})\cap\G _\act \zero = \Theta(s,\X \act s).
  $$

\Proof
  Given $x$ in $\Dm \act {s^*s}$, notice that by \cite [Proposition 4.11]{actions}, one has that $[s,x]$ lies in the
unit space $\G _\act \zero $ if and only if $[s,x]=[e,x]$, for some idempotent $e$, with $x\in\Dm \act e$.  In this case,
there exists an idempotent $f$, such that $x\in\Dm \act f$, and $sf=ef$, so we see that $ef\leq s$, and
  $$
  x\in\Dm \act {e}\cap \Dm \act {f} = \Dm \act {ef} \subseteq \X \act s,
  $$
  whence $[s,x]\in \Theta(s,\X \act s)$.  This proves that $\Theta(s,\Dm \act {s^*s})\cap\G _\act \zero \subseteq \Theta(s,\X \act s)$.

Conversely,
if $x\in\X \act s$, we may choose an idempotent $e$ with $x\in \Dm \act e$, and $se=e$.  It follows that
  $$
  [s,x]=[e,x] \in \Theta(s,\Dm \act {s^*s})\cap\G _\act \zero ,
  $$
  thus proving the reverse inclusion.  \endProof

We are indebted to Benjamin Steinberg for an interesting discussion which helped shape the following characterization of
Hausdorffness for the groupoid of germs.

\state Theorem \label ThmHausdorffAction
  Let $\act $ be the action of an inverse semigroup $\S $ on a locally compact Hausdorff space $X$, satisfying the
conditions of \ref{DefineAction}.  Also let $\G _\act $ be the associated groupoid of germs.  Then, the following
are equivalent:
  \initem
  \nitem $\G _\act $ is Hausdorff,
  \nitem for every $s\in\S $, one has that $\X \act s$ is closed relative to $\Dm \act {s^*s}$.

\Proof Assuming that $\G _\act $ is Hausdorff, we have by \ref{ThmHausdorffGroupoid} that $\G _\act \zero $ is
closed, whence $\Theta(s,\X \act s)$ is closed in $\Theta(s,\Dm \act {s^*s})$ by \ref{TrivialGermsInSlice}.  Observing that
the source map
  $$
  \src : [s,x]\in\Theta(s,\Dm \act {s^*s})\mapsto x\in\Dm \act {s^*s}
  $$
  is a homeomorphism by \cite [Proposition 4.18]{actions}, we deduce that $\src \big(\Theta(s,\X \act s)\big)$ is closed in $\src
(\Theta(s,\Dm \act {s^*s}))$, which is to say that $\X \act s$ is closed in $\Dm \act {s^*s}$.

In order to prove the converse, notice that by the argument used just above, condition (2) implies that $\Theta(s,\X \act s)$
is closed in $\Theta(s,\Dm \act {s^*s})$, so the set
  $$
  \Theta(s,\Dm \act {s^*s})\setminus\Theta(s,\X \act s)
  $$
  is open relative to $\Theta(s,\Dm \act {s^*s})$, and hence also open in $\G _\act $.  Moreover, we have
  $$
  \G _\act \setminus\G _\act \zero =
  \big(\medcup _{s\in\S }\Theta(s,\Dm \act {s^*s})\big)\setminus\G _\act \zero \$=
  \medcup _{s\in\S }\Theta(s,\Dm \act {s^*s})\setminus\G _\act \zero \={TrivialGermsInSlice}
  \medcup _{s\in\S }\Theta(s,\Dm \act {s^*s})\setminus\Theta(s,\X \act s),
  $$
  which is therefore an open set, proving that $\G _\act \zero $ is closed, whence $\G _\act $ is Hausdorff thanks to
\ref{ThmHausdorffGroupoid}.
  \endProof

See also \cite[Theorem 5.17]{SteinbAdv} for the case of actions such that all the $\Dm\alpha e$ are clopen subsets.

Recall from \acite [Theorem 13.3] that, given an inverse semigroup $\S $, the groupoid of germs for the action $\theta:\S
\curvearrowright \Eth $ introduced in \ref{IntroduceTightAction}, is denoted $\Gth $.  We will now use
  \ref {ThmHausdorffAction} to give a characterization of the Hausdorff property for this groupoid.  We
are thankful to Charles Starling for an interesting discussion from where the following result came to life.

\state Theorem \label HausdFiniteCover 
  Let $\S $ be an inverse semigroup (with zero).  Then the following are equivalent:
  \izitem
  \zitem $\Gth $ is Hausdorff,
  \zitem for every $s$ in $\S $, the ideal $\J _s$ defined in \ref{DefineXs} admits a finite cover.

\Proof
  Assuming (ii) we will prove $\Gth $ to be Hausdorff via \ref{ThmHausdorffAction}, so we need to show that the set
$\X \theta s$ introduced in \ref{DefineXs} is closed in $\Dm \theta {s^*s}$, for any $s$ in $\S $.
  Given $s$, let $C$ be a finite cover for $\J _s$.  Then
  $$
  \X \theta s \={XsisDJs}
  \D \theta{\J _s} \={CoversAndCovers} \D \theta C = \medcup _{e\in C}\Dm \theta e.
  $$

Recall that $\Dm \theta e$ is compact for every $e$ in $\E $, as discussed right after \cite [Definition 10.2]{actions}.
Being a finite union of compact sets, $\X \theta s$ is also compact, and hence closed in $\Dm \theta {s^*s}$.  So $\Gth $ is
Hausdorff by \ref{ThmHausdorffAction}.

Conversely, assuming that $\Gth $ is Hausdorff, pick any $s$ in $\S $ and let us produce a finite cover for $\J _s$.
Employing \ref{ThmHausdorffAction} once more, we have that $\X \theta s$ is closed in $\Dm \theta {s^*s}$, and since $\Dm \theta
{s^*s}$ is compact, we have that $\X \theta s$ is itself compact.  On the other hand
  $
  \{D_e\}_{e\in\J _s}
  $
  is clearly an open covering (in the usual topological sense) for $\X \theta s$, so we may extract a finite sub-covering,
say
  $$
  \X \theta s = \medcup _{i=1}^n\Dm \theta {e_i}.
  \equationmark XthetaFiniteUnion
  $$
  Letting
  $$
  C = \{e_1,e_2,\ldots,e_n\},
  $$
  we then have that \ref{XthetaFiniteUnion} translates into
  $$
  \D \theta{\J _s} = \D \theta C.
  $$
  So, using \ref{CoversAndCovers} again, we deduce that $C$ is a finite cover for $\J _s$.
 \endProof

A similar characterization of Hausdorffness for Paterson's universal groupoid of an inverse semigroup was obtained by
Steinberg in \cite{SteinbergPARS}.

Notice that the condition in the above statement that $\J _s$ admit a finite cover is always true if $s$ is idempotent,
since in this case $\{s\}$ is a finite cover for $\J _s$.  Thus it is only relevant to require a finite cover for $\J _s$,
when $s\in\S \setminus \E $.  Given such an $s$, the existence of any nonzero idempotent $e$ in $\J _s$ is a
counter-example to the condition that $\S $ is E*-unitary.  So, requiring a finite cover for $\J _s$, above, may be
interpreted as requiring the E*-unitary property not to fail too badly.  In particular we have:

\state Corollary
  If $\S $ is an E*-unitary inverse semigroup then $\Gth $ is Hausdorff.

\section Topologically free actions

\label TopFreeSect
  In this section we will study fixed points for inverse semigroup actions, paying special attention to situations where
the number of fixed points is small in a certain well defined technical sense, generalizing the well known notion of
\"{topologically free} group actions.  We will then show that our condition is equivalent to the groupoid of germs being
\"{essentially principal}.  Finally we will discuss conditions on an inverse semigroup for the standard action on the
tight spectrum of its idempotent semi-lattice to be topologically free, and hence for the associated tight groupoid to
be essentially principal.

\definition \label DefTopFree
  Given an action $\act :\S \curvearrowright X$, let $s\in\S $, and let $x\in \Dm \act {s^*s}$.
  \initem
  \nitem When $\act _s(x)=x$, we will say that $x$ is a \"{fixed point} for $s$.
  \nitem If there exists $e\in\E $, such that $e\leq s$, and $x\in\Dm \act e$, we will say that $x$ is a \"{trivial fixed point}
for $s$.
  \nitem We say that $\act $ is a \"{free} action, if every fixed point for every $s$ in $\S $ is trivial.
  \nitem We say that $\act $ is a \"{topologically free} action, if for every $s$ in $\S $, the interior of the set of
fixed points for $s$ consists of trivial fixed points.

We will not use the notion of free actions in this work, having presented it above mainly for comparison purposes.

Observe that in case $x$ is a trivial fixed point for $s$, and $e$ is as in \ref{DefTopFree.2} then $se=e$, so
  $$
  \act _s(x) = \act _s\big(\act _e(x)\big) = \act _{se}(x) = \act _e(x) = x,
  $$
  so in particular, $x$ is a fixed point for $s$.  Moreover, since $e\leq s^*s$, by \ref{CompareIniFin}, we have that
$\Dm \act e\subseteq\Dm \act {s^*s}$, and then every $y\in\Dm \act e$ is seen to be a trivial fixed point for $s$.  As a conclusion
we see that the set of trivial fixed points for $s$ is open, and hence it is necessarily contained in the interior of
the set of fixed points for $s$.  In other words, denoting the set of fixed points for $s$ by $F_s$, and the set of
trivial fixed points by $\TF _s$, we always have that
  $$
  \TF _s\subseteq \interior {F_s},
  \equationmark NotationTF
  $$
  while the reverse inclusion holds provided $\act $ is topologically free.

Also notice that if $e$ is an idempotent element in $\S $, then every $x$ in
  $\Dm \act e$
  is a trivial fixed point for $e$.

In the special case of E*-unitary inverse semigroups we have:

\state Proposition
  Given an E*-unitary inverse semigroup $\S $, and
  an action $\act :\S \curvearrowright X$, then
  only idempotent elements may have any trivial fixed points.

\Proof Suppose that an element $s\in\S $ admits  a trivial fixed point $x$.  So there exists an idempotent $e$ such that
$x\in\Dm \act e$, and $se=e$.  It follows that $e\leq s$, and $e\neq0$, because $x\in\Dm \act e$, so we conclude that $s$ is
idempotent.  \endProof

The above notions of freeness and topological freeness therefore become greatly simplified for E*-unitary inverse
semigroups, resembling the corresponding notions for group actions:

\state Proposition \label TheoremRuyTwo
  Given an action $\act :\S \curvearrowright X$, where $S$ is E*-unitary, one has that:
  \initem
  \nitem $\act $ is free if and only if, for every $s\in\S \setminus \E $, the set of fixed points for $s$ is empty.
  \nitem $\act $ is topologically free if and only if, for every $s\in\S \setminus \E $, the set of fixed points for $s$
has empty interior.

\Proof Left for the reader. \endProof

The set $\X \act s$, which played such a crucial role in the study of Hausdorff groupoids (see \ref{ThmHausdorffAction}),
is also relevant regarding trivial fixed points:

\state Proposition \label XsIsIsTrivFixPoints
  Given an action $\act :\S \curvearrowright X$, and given $s$ in $\S $, one has that $\X \act s$ is precisely the set
of trivial fixed points for $s$.

\Proof
  For any $x\in\Dm \act {s^*s}$, one has by definition that $x$ is a trivial fixed point for $s$ if and only if there is
some $e$ in $\E $ such that $e\leq s$, and $x\in\Dm \act e$, and this happens to be precisely the definition of $\X \act s$.
\endProof

Let us now suppose we are given a locally compact \'etale groupoid $\G $,
  with range and source maps $r$ and $d$, respectively, and
  unit space $\G \zero $.  The following are well established notions in the theory of groupoids:

\definition \label {DefEssPrinc}
  (\cite [Definition 3.1]{RenaultCartan})
  \initem
  \nitem The \"{isotropy group bundle} of $\G $ is defined to be the set
  $$
  \G ' = \{\gamma \in\G : d(\gamma) = r(\gamma)\}.
  $$
  \nitem For any $x\in\G \zero $, the \"{isotropy group} of $x$ is defined to be the set
  $$
  \G (x)=\{ \gamma\in\G : d(\gamma)=r(\gamma)=x\}.
  $$
  \nitem $\G $ is said to be \"{principal} if $\G '=\G \zero $.
  \nitem $\G $ is said to be \"{essentially principal} if the interior of $\G '$ coincides with $\G \zero $.

We will not use the notion of principal groupoids in this work, having presented it above mainly for comparison
purposes.

By a result of Renault \cite [Proposition 3.1]{RenaultCartan}, if $\G $ is moreover second countable and Hausdorff,
and $\G \zero $ has the Baire property, then $\G $ is essentially principal if and only if the set consisting of the
points $x$ in $\G \zero $ with trivial isotropy (meaning that $\G (x)=\{x\}$), is dense in $\G \zero $.

In what follows we will relate the important notions defined in \ref{DefTopFree.4} and \ref{DefEssPrinc.4},
via the groupoid of germs.

\state Theorem \label TheoremRuyThree
  Given an action $\act :\S \curvearrowright X$, the corresponding groupoid of germs $\G _\act $ is essentially
principal if and only if $\act $ is topologically free.

  \Proof
  Given $s$ in $\S $, let $x$ be an interior fixed point for $s$.  So there exists an open neighborhood $U$ of $x$,
contained in $\Dm \act {s^*s}$, and consisting of fixed points for $s$.  If $y$ is in $U$, consider the germ $\gamma=[s,y]$.
Then
  $$
  \ran (\gamma) = \act _s(y) = y = \src (\gamma),
  $$
  so $\gamma\in\G _\act '$.  This implies that $\Theta(s,U)$ is contained in $\G _\act '$, and since the former is an open set, we
have that $\Theta(s,U)$ is in fact contained in the interior of $\G _\act '$.  Assuming that $\G _\act $ is essentially
principal, we then deduce that
  $$
  \Theta(s,U)\subseteq\G _\act \zero .
  $$

  In particular $[s,x]$ lies in $\G _\act \zero $, so $[s,x]=[e,x]$, for some idempotent $e\in\E $.  The condition for
equality of germs then gives an idempotent $f$ in $\E $, such that $x\in\Dm \act f$, and $sf = ef$.  Since
  $$
  x\in\Dm \act e\cap\Dm \act f = \Dm \act {ef},
  $$
  and
  $$
  sef = sfe = ef,
  $$
  we conclude that $x$ is a trivial fixed point for $s$.  This shows that every interior fixed point is trivial, and
hence that $\act $ is topologically free.

Conversely, assume that $\act $ is topologically free, and let $\gamma$ lie in the interior of $\G _\act '$.  By definition
of the topology on the groupoid of germs one may then choose $s\in\S $, and an open set $U\subseteq\Dm \act {s^*s}$, such that
  $$
  \gamma\in\Theta (s,U)\subseteq\G _\act '.
  $$
  In particular we may write $\gamma = [s,x]$, for some $x$ in $U$.
  Given any $y$ in $U$, we then have that $[s,y]\in\G _\act '$, so
  $$
  \act _s(y) = \ran ([s,y]) = \src ([s,y]) = y,
  $$
  and we see that $y$ is a fixed point for $s$.  It follows that $U$ is contained in the set of fixed points for $s$.
In particular $x$ is an interior fixed point, and hence, by hypothesis, $x$ is a trivial fixed point.  This means that
there exists an idempotent $e\in\E $, such that $x\in\Dm \act e$, and $se=e$, and consequently
  $$
  \gamma = [s,x] = [e,x] \in \G _\act \zero .
  $$
  This shows that the interior of $\G _\act '$ is contained in $\G _\act \zero $, which is to say that $\G _\act $ is
essentially principal.
  \endProof


The remainder of this section will be devoted to  characterizing  topological freeness for the action of an inverse semigroup
on its tight spectrum.  The next few concepts will  be useful for this purpose.

\definition \label DefineQfixed
  Let $\S $ be an inverse semigroup, and let $s\in\S $.  Given an idempotent $e\in\E $ such that $e\leq s^*s$, we will say that:
  \initem
  \nitem $e$ is \"{fixed} under $s$, if $se=e$,
  \nitem $e$ is \"{weakly-fixed} under $s$, if $sfs^*\Cap f$, for every nonzero idempotent $f\leq e$.

Observe that to say that $e$ is fixed under $s$ is the same as saying that $e\leq s$, as already pointed out  in
\ref{compareIdempWithNonIdemp}.  Moreover notice that the set of all fixed idempotents is precisely $\J _s$.

If $e$ is fixed by $s$, then
  $$
  ses^* = es^* = (se)^* = e,
  $$
  so $e$ is also fixed under \"{conjugation} by $s$, meaning that $ses^* = e$.  Still assuming that $e$ is fixed
under $s$, observe that for every idempotent $f\leq e$, one has that
  $$
  sf = sef = ef = f,
  $$
  whence $f$ is also fixed under $s$.  By the above argument we then have that $sfs^*=f$, so if in addition $f$
is nonzero, then
  $$
  (sfs^*)f = f \neq 0,
  $$
  meaning that $sfs^*\Cap f$.  This holding for any nonzero $f\leq e$, we see that $e$ is weakly-fixed under $s$.  In other
words, every fixed idempotent is weakly-fixed.

\bigskip For the case of the standard action on the space of tight filters, we have:

\state Lemma \label QuasiFixedVsFixedPoints
  Let $\S $ be an inverse semigroup, and let $s\in\S $.  Given an idempotent $e\in\E $ such that $e\leq s^*s$, the following are
equivalent:
  \izitem
  \zitem $e$ is weakly-fixed under $s$,
  \zitem every tight filter $\xi$ in $\Dm \theta e$ is a fixed point for $s$, relative to the standard action $\theta$ of $\S$ on $\Eth $.

\Proof
  Assuming (i), let $\xi$ be an ultra-filter in $\Dm \theta e$.  We then claim that
  $$
  c\Cap d\for c\in\xi \for d\in\theta_s(\xi).
  $$
  Given $c$ and $d$, as above, we have by \ref{DefineActionOfFilters} that $d\geq sbs^*$, for some $b\in\xi$.  So
  $$
  f := bce \in\xi,
  $$
  from where we deduce that $f\neq0$.  Since $f\leq e$, we have by hypothesis that
  $$
  0 \neq (sfs^*)f \leq (sbs^*)c \leq dc.
  $$

  This proves that $dc\neq0$, as claimed.  Since $\xi$ is supposed to be an ultra-filter, we have by \acite [Lemma 12.3] that
$d\in\xi$, for all $d$ in $\theta_s(\xi)$, which is to say that $\theta_s(\xi)\subseteq\xi$.  Observing that $\theta_s(\xi)$ is also an ultra-filter by
\ref{ActionPreserveUltra}, we deduce that $\theta_s(\xi)=\xi$,  proving that $\xi$ is fixed by $s$.

If $\xi$ is a general tight filter in $\Dm \theta e$, then by \acite [Theorem 12.9] we may write $\xi$ as the limit of a net
$\{\xi_i\}_i$ of ultra-filters.  Observing that $\Dm \theta e$ is open, we may also assume that the $\xi_i$ lie in $\Dm \theta e$.  By
what was said above we then see that the $\xi_i$ are fixed by $s$, and hence so is $\xi$, by continuity of $\theta_s$.  This
proves (ii).

Conversely, let $f$ be a nonzero element with $f\leq e$. Using Zorn's Lemma, let $\xi$ be an ultra-filter containing $f$,
which is therefore also a tight filter by \acite [Proposition 12.7].  Since $f\leq e$, we see that $e\in\xi$, whence $\xi\in\Dm \theta
e$.  Using hypothesis (ii) we then have that $\theta_s(\xi)=\xi$, so it follows that $sfs^*\in\xi$, and then necessarily $sfs^*\Cap f$.
This proves that $e$ is weakly-fixed under $s$.  \endProof

Our main result regarding topological freeness for the standard action on the space of tight filters is in order.  The
reader should observe that condition (iii) below is exactly Definition \ref{DefTopFree.4} of topological freeness,
but for the fact that it only refers to ultra-filters.

\state Theorem \label TightTopFree 
  Given an inverse semigroup $\S $ (with zero), consider the following statements:
  \izitem
  \zitem The standard action $\theta:\S\curvearrowright \Eth $ is topologically free.
  \zitem For every $s$ in $\S $, and for every  $e$ in $\E $ which is weakly-fixed under $s$, there exists a
finite cover for $e$ consisting of fixed idempotents.
  \zitem For every $s$ in $\S $, and for every $\xi\in\Dm\theta {s^*s}$ which is an interior fixed point for $s$, as well as an
\underbar {ultra-filter}, one has that $\xi$ is trivially fixed by $s$.
  \medskip \noindent Then {\rm (i)$\Rightarrow $(ii)$\Rightarrow $(iii)}.  If moreover every tight filter is an
ultra-filter\fn
  {Semi-lattices in which every tight filter is an ultra-filter have been  called \"{compactable} semi-lattices \cite [Theorem
2.5]{LawsonCompactable} and they occur quite often%
  .},
  or if $\S $ satisfies the equivalent conditions of \ref{HausdFiniteCover}, then also {\rm
(iii)$\Rightarrow $(i)}.

\Proof
  (i)$\Rightarrow $(ii):
  Let $s\in\S $ and let $e$ be an idempotent element weakly-fixed under $s$.  We then claim that
  $$
  \Dm \theta e =
  \bigcup _{f\in\J _e\cap\J _s}\Dm \theta f.
  \equationmark DeInDf
  $$

By \ref{QuasiFixedVsFixedPoints} we have that
  $\Dm \theta e\subseteq F_s$,
  where $F_s$ denotes the set of fixed points for $s$.  Noticing that $\Dm \theta e$ is open, it follows that
  $$
  \Dm \theta e \subseteq \interior {F_s}.
  $$

  Assuming that $\theta$ is topologically free, we then have that $\Dm \theta e$ consists of trivial fixed points.
  Therefore for any $\xi$ in $\Dm \theta e$, there exists some $f$ in $\E $, such that $sf=f$, and $\xi\in\Dm \theta f$.
  Given that $\xi$ is also in $\Dm \theta e$, we have that $\xi\in\Dm \theta {ef}$, and clearly $sef = ef$.  Therefore $ef$ lies in $\J
_e\cap\J _s$, proving the inclusion ``$\subseteq$'' in \ref{DeInDf}.  Since the reverse inclusion holds trivially, the claim is
proved.

  Using the fact that $\Dm \theta e$ is compact, we may then find a finite set
  $$
  \{f_1,f_2,\ldots,f_n\} \subseteq \J _e\cap\J _s,
  $$
  such that
  $$
  \Dm \theta e=\medcup _{i=1}^n\Dm \theta {f_i}.
  $$

  We then have by \ref{CoversAndCoversForIdempotents} that $\{f_1,f_2,\ldots,f_n\}$ is a cover for $e$.  Since each $f_i$
is in $\J _s$, we have that $f_i$ is fixed by $s$.  This proves (ii).

\medskip \noindent (ii)$\Rightarrow $(iii):
  Pick any $s$ in $\S $, and let $\xi$ be an ultra-filter which is an interior fixed point for $s$.  Then
  %
  by \ref{NicerNeighborhood} there is some $e$ in $\E$ such that
  $$
  \thickmuskip 12mu
  \xi\in U\big(\{e\},\emptyset \big)\cap\Eth\subseteq F_s.
  $$
  Incidentally, the set intersection occurring above coincides with what we have been calling $\Dm\theta e$ starting with
\ref{IntroduceDomains}, hence
  $$
  \xi\in\Dm\theta e \subseteq F_s.
  $$

  By \ref{QuasiFixedVsFixedPoints} we then have that $e$ is weakly-fixed under $s$, so we may use condition (ii) to
obtain a finite cover $\{f_1,f_2,\ldots,f_n\}$ for $e$ consisting of fixed idempotents.  Again invoking \ref{CoversAndCoversForIdempotents} we get
  $$
  \Dm \theta e=\medcup _{i=1}^n\Dm \theta {f_i},
  $$
  so $\xi\in\Dm \theta {f_i}$, for some $i$, and hence $\xi$ is trivially fixed by $s$.

\medskip \noindent (iii)$\Rightarrow $(i): When all tight filters are ultra-filters, (iii) becomes the very definition of
topological freeness and hence there is nothing to do.  Otherwise, assuming the equivalent conditions of \ref{HausdFiniteCover},
  let $s$ be in $\S $, and let $\xi$ be any element of $\Dm\theta {s^*s} $ which lies in the interior of the fixed point set for $s$.

Using \acite [Theorem 12.9] again, let $\{\xi_i\}_i$ be a net of ultra-filters converging to $\xi$, and we may clearly suppose
that all $\xi_i$'s are also interior fixed points for $s$.

By (iii) we have that the $\xi_i$ are trivial fixed points, which is to say that $\xi_i\in\X \theta s$, as observed in \ref{XsIsIsTrivFixPoints}.
  Given that $\Gth $ is Hausdorff, we have by \ref{ThmHausdorffAction} that $\X \theta s$ is closed in $\Dm \theta{s^*s}$, and hence
  $$
  \xi=\lim _i\xi_i \in\X \theta s,
  $$
  whence $\xi$ is also a trivial fixed point.  This proves that $\theta$ is topologically free.
  \endProof

The following is a useful consequence of our work so far:

\state Corollary \label HausdEssPrin
  Let $\S $ be an inverse semigroup.  Then $\Gth $ is both Hausdorff and essentially principal if and only if the
following two conditions hold for every $s$ in $\S $:
  \izitem
  \zitem there exists a finite cover for $\J _s$,
  \zitem for every idempotent $e$ in $\E $ which is weakly-fixed under $s$, there exists a finite cover for $e$
consisting of fixed idempotents.

\Proof
  If $\Gth $ is Hausdorff then (i) holds by \ref{HausdFiniteCover}.  On the other hand, if $\Gth $ is essentially
principal then the standard action $\theta:\S \curvearrowright \Eth $ is topologically free by \ref{TheoremRuyThree} and
hence (ii) holds by the implication (i)$\Rightarrow $(ii) of \ref{TightTopFree}.

Using \ref{HausdFiniteCover} we see that (i) implies that $\Gth $ is Hausdorff and hence the ``equivalent conditions
of \ref{HausdFiniteCover}'' hold, in which case all of the conditions in \ref{TightTopFree} are equivalent.  So
we deduce from (ii) that $\theta$ is topologically free and hence \ref{TheoremRuyThree} implies that $\Gth $ is
essentially principal.
  \endProof

According to a result due to Steinberg, $\Gth$ is essentially principal (or \"{effective}, to use his  terminology) when
$\S$ is 0-disjunctive \cite{SteinbergPARS}.

It is curious that in the ``only if'' part of the above result, the two conditions ``Hausdorff'' and ``essentially principal''
separately imply conditions (i) and (ii), respectively.  However, in proving the converse, apparently  (i) and (ii)
are both needed to get ``essentially principal''.

\section Irreducible actions and minimality

In this short section we will study irreducibility of inverse semigroup actions versus minimality of the corresponding
groupoid of germs.  The main result of this section, namely \ref{MinimalVsIrredTight}, is a characterization of
minimality for the tight groupoid associated to a general inverse semigroup (with zero).  As before, we will fix an
inverse semigroup action $\act :\S \curvearrowright X$ satisfying the conditions laid down in \ref{DefineAction}.

\definition
  \initem
  \nitem Given $x$ and $y$ in $X$, we say that $x$ and $y$ are \"{trajectory-equivalent} under $\act $, in symbols,
$x\sim_{\act }y$, if there exist $s\in\S $ such that $x\in\Dm \act {s^*s}$ and $\act _s(x)=y$.
  \nitem We say that a subset $W$ of $X$ is \"{invariant} under $\act $ if, for every $w\in W$ and $x\in X$, one has that
$w\sim_{\act }x$ implies that $x\in W$.
  \nitem We say that $\act $ is \"{irreducible} if there are no  open invariant subsets of $X$, other than the empty set
and $X$, itself.

It is easy to see that trajectory-equivalence is an equivalence relation.  Also, given a subset $W$ of $X$, it is
elementary to check that $W$ is invariant if and only if
  $$
  \act _s(W\cap\Dm \act {s^*s})\subseteq W \for s\in\S .
  $$

Given any subset $V\subseteq X$, consider the set
  $$
  \orb V = \medcup _{s\in\S } \act _s(V\cap\Dm \act {s^*s}).
  $$
  It is then evident that $\orb V$ is an invariant subset of $X$.  If $V$ is moreover open, then $\orb V$ is clearly also
open.

\state Proposition \label OrbitOfOpen
  A necessary and sufficient condition for $\act $ to be irreducible is that $\orb V=X$, for every nonempty open subset
$V\subseteq X$.

\Proof Left for the reader. \endProof

The corresponding, and  well established notion of minimality for groupoids is as follows:

\definition 
  Given a groupoid $\G $ with range and source maps $r,d:\G \rightarrow \G ^{(0)}$, we say that a subset $U$ of $\G
^{(0)}$ is \"{invariant} if, for every $\gamma$ in $\G $, one has that
  $$
  d(\gamma)\in U \iff r(\gamma)\in U.
  $$
  A groupoid $\G $ is said to be \"{minimal} if the only invariant open subsets of $\G ^{(0)}$ are the empty set and $\G
^{(0)}$ itself.

\state Proposition \label MinimalVsIrred
  Given an action $\act :\S \curvearrowright X$, let $\G _\act $ be the corresponding groupoid of germs.  Then,
identifying $\G ^{(0)}$ with $X$, as usual, we have that the above two notions of invariance for subsets of $X$ agree with
each other.  In particular $\act $ is irreducible if and only if $\G _\act $ is minimal.

\Proof Left for the reader.  \endProof

For the case of the standard action on the space of tight filters, one has:

\state Theorem \label MinimalVsIrredTight 
  Let $\S $ be an inverse semigroup (with zero).  Then the following are equivalent
  \izitem
  \zitem The standard action $\theta:\S \curvearrowright \Eth $ is irreducible,
  \zitem $\Gth $ is minimal,
  \zitem for every nonzero $e$ and $f$ in $\E $, there are $s_1,s_2,\ldots,s_n$ in $\S $, such that $\{s_ifs_i^*\}_{1\leq i\leq n}$ is an
outer cover for $e$.

\Proof
  The equivalence between (i) and (ii) of course follows from \ref{MinimalVsIrred}.

\medskip \noindent (i)$\Rightarrow $(iii): Assuming (i), and given $e$ and $f$ as in (iii), observe that
  the orbit of $\Dm \act f$ coincides with $\Eth $ by \ref{OrbitOfOpen}.  In particular we have that
  $$
  \Dm \act e \subseteq \orb {\Dm \act f} =
  \medcup _{s\in\S } \act _s(\Dm \act f \cap\Dm \act {s^*s}).
  $$

The sets appearing in the right-hand-side above have a nicer description as follows:
  $$
  \act _s(\Dm \act f \cap\Dm \act {s^*s}) =
  \act _s(\Dm \act {fs^*s}) =
  \Dm \act {s(fs^*s)s^*} =
  \Dm \act {sfs^*}.
  \equationmark ActionOnIntersDom
  $$
  Therefore $\{\Dm \act {sfs^*}\}_{s\in\S }$ is an open cover for the compact set $\Dm \act e$, and hence there is a
finite sub-cover, say
  $$
  \Dm \act e \subseteq \medcup _{i=1}^n \Dm \act {s_ifs_i^*}.
  $$

  We then conclude from \ref{CoversAndCoversForIdempotents} that  $\{s_ifs_i^*\}_{1\leq i\leq n}$ is an outer cover for $e$, hence
proving point (iii).

\medskip \noindent (iii)$\Rightarrow $(i): Given a nonempty open invariant subset $U\subseteq\Eth $, our task is to show that
necessarily $U=\Eth $.  By \acite [Theorem 12.9] there exists an ultra-filter $\xi$ in $U$, and by \ref{NicerNeighborhood}
there is a (necessarily nonzero) idempotent $f$ such that
  $$
  \xi\in\Dm \act f\subseteq U.
  $$

  It then follows that $\orb {\Dm \act f}\subseteq U$ and, in order to complete the proof, it is enough to prove that
  $$
  \Eth \subseteq\orb {\Dm \act f}.
  $$

  Given $\eta\in\Eth $, pick any $e$ in $\eta$, and use (iii) to obtain $s_1,s_2,\ldots,s_n$ in $\S $, such that $\{s_ifs_i^*\}_{1\leq i\leq n}$
is a cover for $e$.  Therefore
  $$\thickmuskip =14mu
  \eta \in \Dm \act e\explain\subseteq{CoversAndCoversForIdempotents}
  \medcup _{i=1}^n \Dm \act {s_ifs_i^*} \={ActionOnIntersDom}
  \medcup _{i=1}^n \act _{s_i}(\Dm \act f \cap\Dm \act {s_i^*s_i}) \subseteq
  \orb {\Dm \act f}.
  $$
  This concludes the proof.  \endProof

Another sufficient condition for the minimality of $\Gth$, found by Steinberg, is that $\S$ be 0-simple \cite{SteinbergPARS}.

\section Local Contractiveness

\label LocContrSection
  Besides the properties of being Hausdorff, essentially principal, and minimal, a property of \'etale groupoids which has
attracted a fair amount of interest is local contractiveness, since these together imply that the reduced groupoid
C*-algebra is simple and purely infinite by \cite[Theorem 4.4]{nhausdorff} and \cite[Proposition 2.4]{AdelaR} (see also
\cite{SimpleGroupoid}).  With this as our main motivation, we will now discuss local contractiveness, both at the level
of groupoids and of inverse semigroup actions.

Corollary \ref{CharacGTightContractive}, the main result in this section, is a characterization of local
contractiveness for the tight groupoid of an inverse semigroup.  Unfortunately it  has a little hitch,
since it relies on the assumption that every tight filter is an ultra-filter, a problem already encountered in
\ref{TightTopFree}.

We begin by recalling Anantharaman-Delaroche's definition of the main concept used in  the present section.

\definition \label LocContrGpd
  \cite[2.1]{AdelaR} \
  Let $\G$ be a locally compact \'etale groupoid.  We say that $\G$ is \"{locally contracting} if, for every nonempty open subset
$U\subseteq\G\zero$, there exists an open subset $V\subseteq U$ and an open bissection $S\subseteq\G$, such that $\clos V\subseteq S\inv S$, and
  $
  S\clos V S\inv\subsetneqq V.
  $

The appropriate definition for inverse semigroup actions seems to be the following:

\definition \label LocContrAct
  We will say that an action $\act :\S \curvearrowright X$ is \"{locally contracting} if, for every nonempty open subset
$U\subseteq X$, there exists an open subset $V\subseteq U$ and an element $s$ in $\S$, such that
  $\clos V\subseteq\Dm\act{s^*s}$, and
  $\act_s(\clos V)\subsetneqq V$.

\state Proposition \label LocContrForActionsAndGroupoid
  Given a locally contracting action $\act :\S \curvearrowright X$, one has that the corresponding groupoid of germs
$\G_\act$ is locally contracting.

\Proof
Let $U$ be a nonempty open subset of the unit space of $\G_\act$, which we identify with $X$ via \ref{UnitSpace}, as usual.    Assuming that
$\act$ is locally contracting, let
$V$ and $s$ be as in \ref{LocContrAct}.  Considering the bissection
  $$
  S := \Theta(s,\Dm\act {s^*s}),
  $$
  observe that
  $$
  S\inv S = \Dm\act {s^*s} \supseteq \clos V,
  $$
  and
  $$
  S\clos VS\inv = \act_s(\clos V) \subsetneqq V,
  $$
  proving that $\G_\act$ is locally contracting.
  \endProof

We have not been able to determine if the converse of the above result also holds.  The fact that the groupoid of germs
has a lot more bissections than simply the $\Theta(s,\Dm\act {s^*s})$ makes one wonder whether local contractivity for
$\G_\act$ should imply the same for $\act$.  Nevertheless, in our future application of
\ref{LocContrForActionsAndGroupoid}, we will fortunately not need the converse implication.

\definition \label LocContrSGroup
  We will say that an inverse semigroup $\S$ is \"{locally contracting} if, for every nonzero idempotent $e$ in $\S$,
there exists $s$ in $\S$, and a finite subset
  $$
  F:= \{f_0,f_1,\ldots,f_n\} \subseteq \E,
  $$
  with $n\geq0$, such that, for all $i=0,\ldots,n$, one has
  \izitem
  \zitem $0\neq f_i\leq e s^*s$,
  \zitem $F$ is an outer cover for $sf_is^*$,
  \zitem $f_0sf_i=0$.

The best relation we found  between these  homonymic conditions is as follows:

\state Theorem \label LocalContrEquiv
  Given an inverse semigroup $\S$, consider the following
statements:
  \iaitem
  \aitem $\S$ is a locally contracting inverse semigroup,
  \aitem the standard action of\/ $\S$ on $\Eth$ is  locally contracting.
  \iskip Then {\rm (a)$\Rightarrow$(b)}.  If moreover every tight filter on $\E$ is an ultra-filter, then also {\rm
(b)$\Rightarrow$(a)}.

\Proof 
  Suppose $\S$ is locally contracting.  Given a nonempty subset $U\subseteq\Eth$, choose an ultra-filter $\xi$ in $U$.  Employing
\ref{NicerNeighborhood} we then see that  there is an idempotent $e$ in $\E$ such that
  %
  $$
  \xi\in\Dm\theta e\subseteq U.
  $$

  Using  hypothesis (a), let $s$ and    $\{f_0,f_1,\ldots,f_n\}$ be as in \ref{LocContrSGroup}, and define
  $$
  V = \medcup_{i=0}^n\Dm\theta {f_i}.
  $$

  Since each $f_i\leq e$, it is clear that $V\subseteq\Dm\theta e\subseteq U$.  Moreover each $f_i\leq s^*s$, and $V$ is clopen, so
  $$
  \clos V= V\subseteq\Dm\theta {s^*s}.
  $$

  In order to complete the proof it is now enough to prove that $\theta_s(\clos V)\subsetneqq V$.  With this purpose in mind notice
that for each idempotent $f$ in $\E$, one has that $\theta_s\circ\theta_f = \theta_{sf}$, so
  $
  \theta_s(\Dm\theta f)
  $
  coincides with the range of $\theta_{sf}$, also known as $\Dm\theta {sfs^*}$.  Therefore,
  $$
  \theta_s(V) =
  \medcup_{i=0}^n\theta_s(\Dm\theta {f_i}) =
  \medcup_{i=0}^n\Dm\theta {sf_is^*}.
  $$
  By \ref{LocContrSGroup.ii} and
\ref{CoversAndCoversForIdempotents} we have that each
  $
  \Dm\theta {sf_is^*}
  $
  is contained in $V$, hence
  $$
  \theta_s(\clos V) =   \theta_s(V) \subseteq V.
  $$

  We still need to show that   $\theta_s(V)$ is a proper subset of $V$, and to do so we recall from
\ref{LocContrSGroup.iii} that $f_0sf_i=0$, for all $i$.  Consequently
  $f_0\perp sf_is^*$, which implies that $\Dm\theta {f_0}$ is disjoint from  $\Dm\theta {sf_is^*}$.  So
  $$
  \Dm\theta {f_0} \subseteq
  V \setminus \big(\medcup_{i=0}^n\Dm\theta {sf_is^*}\big) =
  V \setminus \theta_s(V) =
  V \setminus \theta_s(\clos V),
  $$
  and since   $\Dm\theta {f_0}$ is nonempty (by Zorn's Lemma it contains an ultra-filter) we see  that
$\theta_s(\clos V)\subsetneqq V$.  This proves   that $\theta$ is a locally contracting action.

Conversely, assume (b) and also  that every tight filter in $\E$ is an ultra-filter.
Given a nonzero $e$ in $\E$, let $U=\Dm\theta e$, and choose $s$ and $V$ as in
\ref{LocContrAct} so that, among other things,
  $$
  \theta_s(\clos V) \subsetneqq V\subseteq U.
  $$

  For each $\xi$ in $\theta_s(\clos V)$, choose a neighborhood of $\xi$ contained in $V$.  By hypothesis we have that $\xi$ is an
ultra-filter, and by \ref{NicerNeighborhood} we may suppose that such a neighborhood  is of the form
$\Dm\theta {f_\xi}$, for some $f_\xi$ in $\E$, so
  $$
  \xi\in\Dm\theta {f_\xi}\subseteq V.
  $$
  The $\Dm\theta {f_\xi}$ then evidently form an open cover for $\theta_s(\clos V)$.
  Since
  $$
  V\subseteq U\cap\Dm\theta {s^*s} =\Dm\theta e\cap\Dm\theta {s^*s} = \Dm\theta {es^*s},
  \equationmark VInEssStar
  $$
  we may replace  each $f_\xi$  by
  $$
  f_\xi':= es^*sf_\xi,
  $$
  and hence we may assume that $f_\xi\leq es^*s$.

  Being a closed subset of $\Dm\theta {s^*s}$, observe that $\clos V$ is   compact, and hence so  is $\theta_s(\clos
V)$.  We may then take a finite subcover of the above cover, say
  $$\thickmuskip 12mu
  \theta_s(\clos V)\subseteq\medcup_{f\in F'} \Dm\theta {f} \subseteq V,
  \equationmark LocalContrEquivClosV
  $$
  where $F'$ is a finite set consisting of some of the $f_\xi$.

We next claim that there exists  a nonzero idempotent $f_0\leq es^*s$, such that
  $$
  \Dm\theta {f_0}\subseteq V\setminus \theta_s(\clos V).
  \equationmark LocalContrEquivFZero
  $$

  To see this, first observe that
  $
  V\setminus \theta_s(\clos V)
  $
  is  open and nonempty.
Even  without assuming that all tight filters are ultra-filters, we may use the density of the set formed by the latter
to find some ultra-filter $\xi$ in   $V\setminus \theta_s(\clos V)$.  An application of \ref{NicerNeighborhood}  then
provides  $f_0$ in $\E$ such that
  $$
  \xi \subseteq \Dm\theta {f_0} \subseteq V\setminus \theta_s(\clos V),
  \equationmark FZeroAppears
  $$
  and, again by  \ref{VInEssStar}, we may assume that $f_0\leq es^*s$.  Adding $f_0$  to  $F'$, we form the set
  $$
  F := \{f_0\}\cup F'=\{f_0,f_1,\ldots,f_n\},
  $$
  which is then seen to satisfy \ref{LocContrSGroup.i}.
  Defining
  $$
  W=\medcup_{f\in F} \Dm\theta {f},
  $$
  observe that $W$ is clopen, and that
  $W\subseteq V\subseteq\Dm\theta{s^*s}$.  We moreover have  that
  $$
  \thickmuskip 12mu
  \theta_s(W) \subseteq \theta_s(V) \explain\subseteq{LocalContrEquivClosV} W.
  $$

  In particular, for each $i=0,1,\ldots,n$, we have
  $$
  \thickmuskip 12mu
  \Dm\theta {sf_is^*} = \theta_s(\Dm\theta {f_i}) \subseteq \theta_s(W) \subseteq W = \medcup_{f\in F} \Dm\theta {f},
  $$
  so $F$ is an outer cover for $sf_is^*$ by \ref{CoversAndCoversForIdempotents}, proving \ref{LocContrSGroup.ii}.
As seen above, for all $i$, we have
  $$
  \Dm\theta {sf_is^*} \subseteq \theta_s(W) \subseteq \theta_s(V) \subseteq \theta_s(\clos V),
  $$
  so $\Dm\theta {sf_is^*}$ and $\Dm\theta {f_0}$ are disjoint by \ref{FZeroAppears}, and hence
  $$
  \emptyset = \Dm\theta {f_0} \cap\Dm\theta {sf_is^*} = \Dm\theta {f_0sf_is^*},
  $$
  from where we conclude that $f_0sf_is^*=0$, and hence also
  $$
  f_0sf_i =   f_0ss^*sf_i = f_0sf_i s^*s = 0.
  $$
  This  proves \ref{LocContrSGroup.iii}, and hence we are done.
\endProof

An immediate consequence is as follows:

\state Corollary \label CharacGTightContractive 
  For every  locally contracting inverse semigroup $\S$, one has that  $\Gth$ is a  locally contracting groupoid.

\Proof
  Follows immediately from  \ref{LocalContrEquiv} and \ref{LocContrForActionsAndGroupoid}.
\endProof

As we have already pointed out, both \ref{LocalContrEquiv} and \ref{LocContrForActionsAndGroupoid} do not seem to admit
converses.  Therefore a  converse for \ref{CharacGTightContractive}   seems to be even less likely.

   In  applications of \ref{LocalContrEquiv} it is sometimes possible to verify a set
of conditions which is both stronger and formally easier to check than \ref{LocalContrEquiv.i--iii}, as follows:

\state Proposition \label EasierLocConr
  Suppose that, for every nonzero $e\in\E$, there exists $s\in\S$, and idempotent elements $f_0$ and $f_1$, such that
  \izitem
  \zitem $0 \neq f_0 \leq f_1 \leq e s^*s$,
  \zitem $sf_1s^*\leq f_1$,
  \zitem $f_0sf_1=0$.
  \iskip Then $\S$ is locally contracting.

\Proof
  It is enough to observe that conditions \ref{LocContrSGroup.i--iii} are  satisfied for $F=\{f_0,f_1\}$, as the reader may
easily verify.
  \endProof

As a final result, let us put all of the above to work, in order to obtain a class  of simple, purely infinite
C*-algebras:

\def\redalg{C^*\red\big(\Gth\big)}

\state Theorem
  Let $\S$ be a countable inverse semigroup (with zero) and let us consider the following conditions:
  \iaitem
  \aitem For every $s$ in $\S $, the ideal $\J _s$ defined in \ref{DefineXs} admits a finite cover,  
  \aitem For every $s$ in $\S $, and for every  $e$ in $\E $ which is weakly-fixed under $s$, there exists a
finite cover for $e$ consisting of fixed idempotents.
  \aitem For every nonzero $e$ and $f$ in $\E $, there are $s_1,s_2,\ldots,s_n$ in $\S $, such that $\{s_ifs_i^*\}_{1\leq i\leq n}$ is an
outer cover for $e$.
  \aitem $\S$ is locally contracting.
  \iskip
  Then:
  \izitem
  \zitem Conditions (a+b) imply that every nonzero ideal
  $$
  J\trianglelefteq \redalg
  $$
  has a nonzero  intersection  with the algebra of continuous functions vanishing at infinity on the unit space of $\Gth$.
  \zitem Conditions (a+b+c) imply that   $\redalg$ is a simple C*-algebra.
  \zitem Conditions (a+b+c+d) imply that   $\redalg$ is a purely infinite simple C*-algebra.

\Proof
  All we now have to do is put the pieces together.  By  \ref{HausdEssPrin} we have that $\Gth$ is Hausdorff and
essentially principal, and it is also easy to see that it is second countable as a consequence of the assumption that
$\S$ is countable.   We may then  invoke \cite[Theorem 4.4]{nhausdorff} to obtain (i).

Adding (c) to the above,  we have  by \ref{MinimalVsIrredTight} that $\Gth$ is minimal, so a now standard argument
(see e.g.~\cite[Corollary 29.8]{book}) easily implies that $\redalg$ is simple.

If, on top of it all, we assume (d), then $\Gth$ is contracting by \ref{CharacGTightContractive}, so the result follows
from \cite[Proposition 2.4]{AdelaR}.
\endProof

\references

\Article AdelaR
  C. Anantharaman-Delaroche;
  Purely infinite $C^*$-algebras arising form dynamical systems;
  Bull. Soc. Math. France, 125 (1997), no. 2, 199-225

\Article SimpleGroupoid
  J. Brown, L. O. Clark, C. Farthing and A. Sims;
  Simplicity of algebras associated to \'etale groupoids;
  Semigroup Forum, 88 (2014), 433-452

\Article actions
  R. Exel;
  Inverse semigroups and combinatorial C*-algebras;
  Bull. Braz. Math. Soc., 39 (2008), no. 2, 191-313

\Article nhausdorff
    R. Exel;
    Non-Hausdorff \'etale groupoids;
    Proc. Amer. Math. Soc., 139 (2011), no. 3, 897-907

\Bibitem book
  R. Exel;
  Partial Dynamical Systems, Fell Bundles and Applications;
  {Licensed under a Creative Commons Attribution-ShareAlike 4.0 International License, 351pp, 2014.  Available online from\break
mtm.ufsc.br/$\sim$exel/publications. \rm PDF file md5sum: bc4cbce3debdb584ca226176b9b76924}

\Bibitem EP
  R. Exel and E. Pardo;
  Representing Kirchberg algebras as inverse semigroup crossed products;
  arXiv:1303.6268 [math.OA], 2013

\Bibitem EPTwo
  R. Exel and E. Pardo;
  Graphs, groups and self-similarity;
  arXiv:1307.1120 [math.OA], 2013

\Bibitem EPThree
  R. Exel and E. Pardo;
  Self-similar graphs, a unified treatment of Katsura and Nekrashevych C*-algebras;
  in preparation

\Article KatsuraOne
  T. Katsura;
  A construction of actions on Kirchberg algebras which induce given actions on their $K$-groups;
  J. reine angew. Math., 617 (2008), 27-65

\Bibitem Lawson
  M. V. Lawson;
  Inverse semigroups, the theory of partial symmetries;
  World Scientific, 1998

\Article LawsonCompactable
  M. V. Lawson;
  Compactable semilattices;
  Semigroup Forum, 81 (2010), no. 1, 187-199

\Article NekraJO
  V. Nekrashevych;
  Cuntz-Pimsner algebras of group actions;
  J. Operator Theory, 52 (2004), 223-249

\Bibitem Nmsn
  V. Nekrashevych;
  Self-similar groups;
  Mathematical Surveys and Monographs, \bf 117, \rm  Amer. Math. Soc., Providence, RI, 2005

\Article NC
  V. Nekrashevych;
  C*-algebras and self-similar groups;
  J. reine angew. Math., 630 (2009), 59-123

\Bibitem pat
  A. L. T. Paterson;
  Groupoids, inverse semigroups, and their operator algebras;
  Birkh\umlaut auser, 1999

\Article RenaultCartan
  J. Renault;
  Cartan subalgebras in $C^*$-algebras;
  Irish Math. Soc. Bull., 61 (2008), 29-63

\Article SteinbAdv
  B. Steinberg;
  A groupoid approach to discrete inverse semigroup algebras;
  Adv. Math., 223 (2010), 689-727

\Bibitem SteinbergPARS
  B. Steinberg;
  Associative algebras associated to \'etale groupoids and inverse semigroups;
  PARS -- Partial Actions and Representations Symposium, Gramado,
Brazil, 2014.  Available from \break http://mtm.ufsc.br/~exel/PARS/proceedings.html

\Bibitem SteinbergPreprint
  B. Steinberg;
  Simplicity, primitivity and semiprimitivity of \'etale groupoid algebras with applications to inverse semigroup
algebras;
  arXiv:1408.6014 [math.RA]

\endgroup

\bigskip \noindent \tensc
  Departamento de Matem\'atica;
  Universidade Federal de Santa Catarina;
  88010-970 Florian\'opolis SC;
  Brazil
  \hfill \break \tt (exel@mtm.ufsc.br)

\bigskip \noindent \tensc
  Departamento de Matem\'aticas, Facultad de Ciencias;
  Universidad de C\'adiz, Campus de Puerto Real;
  11510 Puerto Real (C\'adiz);
  Spain
  \hfill \break \tt (enrique.pardo@uca.es)

\close

\bye